\documentclass[11pt]{article}

\usepackage[a4paper,margin=1in]{geometry}

\usepackage{amsmath}
\usepackage{amssymb}
\usepackage{amsfonts}
\usepackage{amsthm}

\theoremstyle{remark}
\newtheorem{remark}{Remark}

\usepackage{graphicx}
\usepackage{algorithm}
\usepackage{algpseudocode}

\usepackage{booktabs}
\usepackage{hyperref}
\usepackage{enumitem}
\usepackage{tikz}
\usepackage{pgfplots}

\pgfplotsset{compat=1.18}

\title{
Deep Learning for Singular PDEs: A Weighted Neural Network Approach
}

\author{
Badr Oulgiht\\
}

\begin{document}

\maketitle


\begin{abstract}

We propose a deep learning framework for the numerical approximation of weakly singular semilinear elliptic equations subject to homogeneous Dirichlet boundary conditions and positivity constraints. These problems are challenging because the nonlinear source term becomes unbounded as the positive solution approaches zero near the boundary.

The proposed approach combines three components. First, the Dirichlet boundary conditions are embedded directly into the neural representation through a hard-constraint formulation. Second, a Softplus transformation preserves positivity in the interior of the computational domain, ensuring that the singular nonlinear term remains well-defined during training. Third, a singularity-aware weighted residual loss emphasizes regions where the predicted solution becomes small and the singular behavior is more pronounced.

The method is evaluated in the weak singularity regime through numerical reference solutions and controlled verification tests. Comparisons with the standard residual formulation are used to assess approximation accuracy and convergence behavior. The numerical experiments illustrate the potential of combining hard boundary enforcement, positivity preservation, and singularity-aware residual weighting for weakly singular elliptic problems.
\end{abstract}

\vspace{0.3cm}

\noindent
\textbf{Keywords:}
Physics-Informed Neural Networks,
Deep Learning,
Singular Elliptic Equations,
Singularity-Aware Learning,
Scientific Machine Learning.


\section{Introduction}

Nonlinear partial differential equations (PDEs) involving singular
nonlinearities arise in a variety of mathematical models and scientific
applications. These problems are characterized by nonlinear source terms
that become unbounded as the solution approaches zero, giving rise to
additional analytical and numerical difficulties compared with equations
involving bounded nonlinearities. In particular, the singular behavior may
affect the regularity of the solution and requires special care when
constructing stable and accurate numerical approximations.

In this work, we consider the singular semilinear elliptic problem

\begin{equation}
\label{eq:intro}
\begin{cases}
-\Delta u=\dfrac{1}{u^{\alpha}},
& x\in\Omega,\\[1mm]
u=0,
& x\in\partial\Omega,\\[1mm]
u>0,
& x\in\Omega,
\end{cases}
\end{equation}

where $\Omega\subset\mathbb{R}^{N}$ is a bounded domain and the singularity
parameter satisfies

\begin{equation}
0<\alpha<1.
\label{eq:intro_alpha}
\end{equation}

This range corresponds to the weak singularity regime. Although the solution
is required to remain strictly positive in the interior of the domain, the
homogeneous Dirichlet condition forces it to approach zero at the boundary.
Consequently, the nonlinear source term

\begin{equation}
u^{-\alpha}
\end{equation}

becomes unbounded as the boundary is approached. This interaction between
boundary vanishing and the singular nonlinear source constitutes the main
difficulty addressed in the present work.

Singular semilinear elliptic equations and related singular boundary-value
problems have been studied extensively in nonlinear analysis. They are
relevant to mathematical models in which the state variable is constrained
to remain positive while the governing equation contains a nonlinear
response that becomes singular as the state approaches zero. From both
analytical and computational perspectives, such equations require numerical
methods capable of handling the singular source term while preserving the
positivity of the approximation.

Classical numerical approaches, including finite difference and finite
element methods, can be employed for singular elliptic equations. However,
the unbounded nonlinear source near regions where the solution approaches
zero may require additional numerical treatment, such as regularization,
mesh refinement, or appropriately designed nonlinear iterative procedures.
The simultaneous treatment of the boundary conditions, positivity, and
singular nonlinearity therefore remains an important computational issue.

In recent years, deep learning methods have emerged as an alternative
framework for the numerical approximation of differential equations. Among
these approaches, Physics-Informed Neural Networks incorporate the governing
differential equation directly into the training objective. Spatial
derivatives of the neural approximation are computed through automatic
differentiation, allowing the differential residual to be evaluated at
collocation points without requiring labeled solution data during training.

Many standard PINN formulations are primarily developed and benchmarked on
problems with sufficiently regular solutions and bounded nonlinearities.
Their application to singular elliptic equations introduces additional
difficulties because the nonlinear source term becomes unbounded as the
solution approaches zero. Moreover, an unconstrained neural network may
produce nonpositive values during optimization, which is particularly
problematic for nonlinear terms of the form $u^{-\alpha}$. Preserving
positivity is therefore not merely a qualitative requirement but an
important component of the numerical formulation.

In this work, we develop a PINN framework specifically adapted to weakly
singular semilinear elliptic equations. The proposed approach combines hard
enforcement of the Dirichlet boundary conditions with a
positivity-preserving neural representation. In addition, a
singularity-aware weighted residual loss is introduced to increase the
contribution of collocation points where the predicted solution becomes
small and the singular nonlinear term has a stronger influence.

The main contributions of this work are summarized as follows:

\begin{itemize}

    \item We develop a neural framework for weakly singular semilinear
    elliptic equations in the regime $0<\alpha<1$, where the nonlinear
    source term becomes unbounded as the solution approaches zero.

    \item We construct a hard-constrained, positivity-preserving neural
    representation that satisfies the Dirichlet boundary conditions exactly
    while maintaining strictly positive predictions in the interior of the
    domain.

    \item We introduce a singularity-aware weighted residual loss that
    increases the contribution of regions where the predicted solution is
    small and the singular nonlinear term becomes dominant.

    \item We evaluate the proposed formulation through numerical benchmarks
    and controlled verification tests, including comparisons with the
    corresponding standard residual formulation.

\end{itemize}

Several recent works have explored neural-network-based methods for the
numerical solution of nonlinear and elliptic partial differential equations,
including Physics-Informed Neural Networks, adaptive loss formulations, and
constrained neural architectures \cite{Sukumar2022,Raissi2019,Sirignano2018,Wang2021Gradient,EYu2018,Yu2022,Wang2024Causal}. These studies
provide useful background for the singularity-aware neural approach developed
in the present work.\\

The remainder of this paper is organized as follows.
Section~2 introduces the mathematical formulation of the weakly singular elliptic problem and discusses the main properties relevant to the numerical method. Section~3 presents the proposed PINN methodology, including the hard-constrained neural representation, positivity preservation, and the singularity-aware weighted loss function. Section~4 reports the numerical experiments and evaluates the proposed method using singular benchmarks, controlled verification tests, and comparisons with the standard residual formulation. Finally, Section~5 presents the conclusions and possible directions for future research.



\section{Mathematical Overview}
\label{sec:mathematical_overview}

We consider the weakly singular semilinear elliptic boundary-value problem

\begin{equation}
\begin{cases}
-\Delta u
=
\dfrac{1}{u^{\alpha}},
& x\in\Omega,
\\[2mm]
u=0,
& x\in\partial\Omega,
\\[2mm]
u>0,
& x\in\Omega,
\end{cases}
\label{eq:model}
\end{equation}

where $\Omega\subset\mathbb{R}^{d}$ is a bounded domain and the singularity
exponent satisfies

\begin{equation}
0<\alpha<1.
\label{eq:alpha_range_math}
\end{equation}

The nonlinear function

\begin{equation}
s\longmapsto s^{-\alpha},
\qquad s>0,
\end{equation}

is singular at the origin. More precisely,

\begin{equation}
\lim_{s\rightarrow0^{+}}
s^{-\alpha}
=
+\infty.
\label{eq:singular_limit}
\end{equation}

\begin{figure}[ht]
    \centering
    \includegraphics[width=0.68\textwidth]{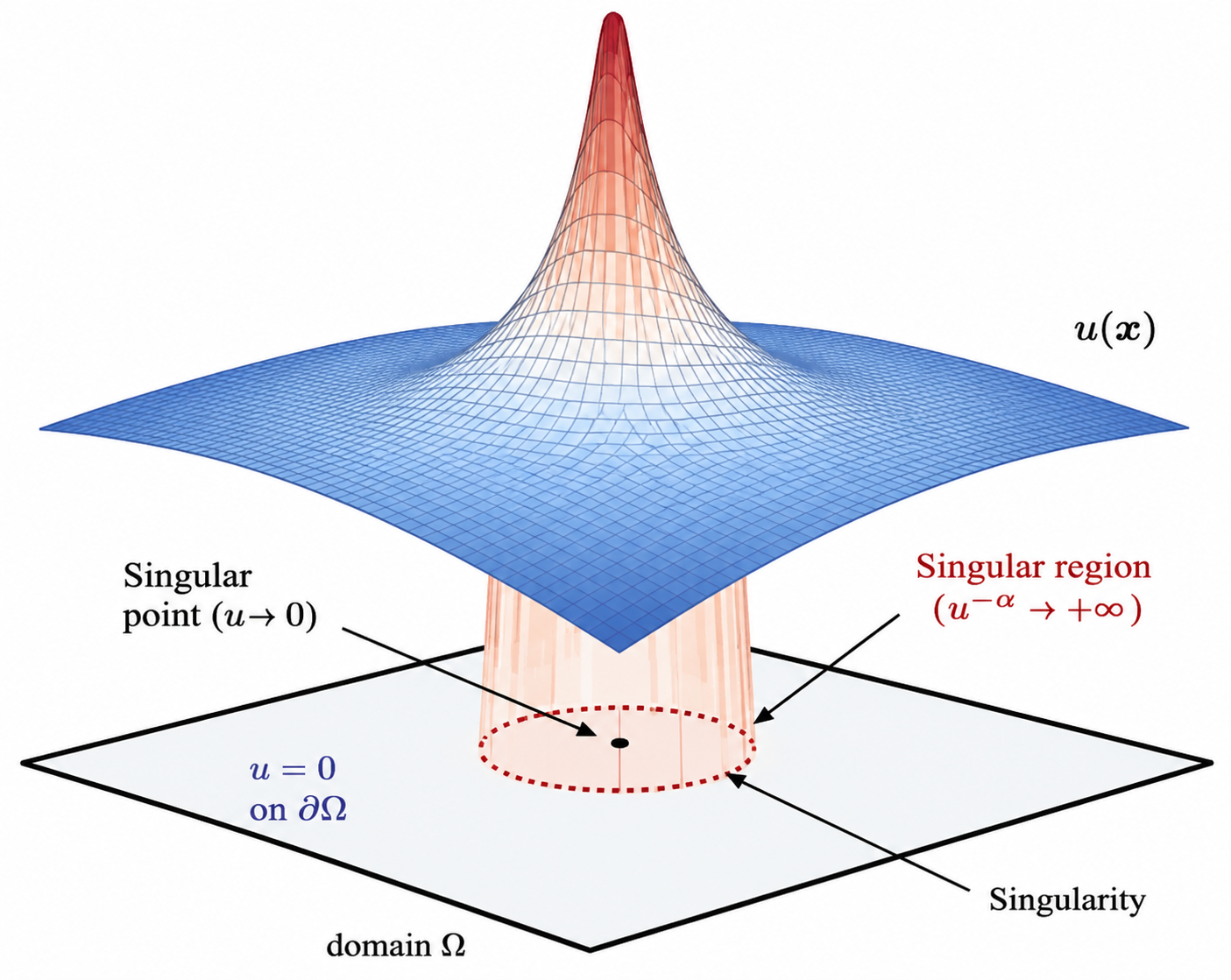}
    \caption{
    Illustration of the singular behavior of the nonlinear function
    $s\mapsto s^{-\alpha}$ for $0<\alpha<1$. The nonlinear term becomes
    unbounded as $s\rightarrow0^{+}$.
    }
    \label{fig:singular_behavior}
\end{figure}

The singular character of the nonlinear source term is illustrated in
Figure~\ref{fig:singular_behavior}. As $s$ approaches zero from the positive
side, the quantity $s^{-\alpha}$ increases without bound. This behavior is
particularly relevant for problem~\eqref{eq:model}, since the homogeneous
Dirichlet condition forces the solution to approach zero at the boundary
while the solution remains positive in the interior of the domain.

Away from the boundary, where the positive solution remains bounded away
from zero, the nonlinear term $u^{-\alpha}$ remains finite. The principal
difficulty arises in regions where the solution becomes small, since the
magnitude of the nonlinear source term increases rapidly as
$u\rightarrow0^{+}$.

From a variational perspective, problem~\eqref{eq:model} is formally
associated with the energy functional

\begin{equation}
\mathcal{J}(u)
=
\frac{1}{2}
\int_{\Omega}
|\nabla u|^{2}\,dx
-
\frac{1}{1-\alpha}
\int_{\Omega}
u^{1-\alpha}\,dx,
\label{eq:energy}
\end{equation}

defined formally for positive admissible functions satisfying the prescribed
boundary condition.

Indeed, the formal first variation of~\eqref{eq:energy} in the direction of
an admissible test function $v$ gives

\begin{equation}
\frac{d}{d\varepsilon}
\mathcal{J}(u+\varepsilon v)
\bigg|_{\varepsilon=0}
=
\int_{\Omega}
\nabla u\cdot\nabla v\,dx
-
\int_{\Omega}
u^{-\alpha}v\,dx.
\end{equation}

After integration by parts, and using the homogeneous boundary condition on
the admissible variations, the corresponding Euler--Lagrange equation is

\begin{equation}
-\Delta u-u^{-\alpha}=0,
\end{equation}

which is equivalent to

\begin{equation}
-\Delta u
=
\frac{1}{u^{\alpha}}.
\end{equation}

The positivity condition plays an essential role in the mathematical and
numerical formulation of the problem. In particular, an approximation that
takes zero or negative values at interior points is incompatible with the
singular nonlinear term $u^{-\alpha}$. A numerical approximation should
therefore account simultaneously for three structural properties:

\begin{equation}
\boxed{
\text{boundary condition}
\quad+\quad
\text{interior positivity}
\quad+\quad
\text{singular residual}.
}
\end{equation}

Let $\hat{u}(x;\theta)$ denote a neural approximation parameterized by
$\theta$. The differential residual associated with
problem~\eqref{eq:model} is

\begin{equation}
R(x;\theta)
=
-\Delta\hat{u}(x;\theta)
-
\frac{1}
{\hat{u}(x;\theta)^{\alpha}}.
\label{eq:mathematical_residual}
\end{equation}

When $\hat{u}(x;\theta)$ becomes small, the magnitude of the nonlinear
contribution

\begin{equation}
\frac{1}
{\hat{u}(x;\theta)^{\alpha}}
\end{equation}

increases. Consequently, the residual may exhibit substantially different
scales across the computational domain. This observation motivates the use
of a residual weighting strategy that assigns additional importance to
regions where the predicted solution is small.

The existence and regularity theory of singular semilinear elliptic
equations provides the mathematical background for the problem considered
here. The objective of the present work, however, is not to establish new
analytical existence or regularity results, but rather to develop and
evaluate a neural numerical framework adapted to the singular structure of
the equation.

The following section introduces the proposed methodology, combining a
hard enforcement of the homogeneous Dirichlet boundary condition, a
positivity-preserving neural representation, automatic differentiation, and
a singularity-aware weighted residual loss.


\section{Methodology}
\label{sec:methodology}

We now introduce the neural framework developed for the numerical
approximation of the weakly singular semilinear elliptic problem
\eqref{eq:model}. The objective is to construct an approximation
$\hat{u}(x;\theta)$ that respects the essential structure of the problem
while allowing the governing differential equation to be incorporated
directly into the training process.

The principal numerical difficulty is associated with the singular nonlinear
term
\[
u^{-\alpha},
\qquad 0<\alpha<1,
\]
which becomes unbounded as the positive solution approaches zero. Since the
homogeneous Dirichlet condition forces the solution to vanish at the
boundary, the neural approximation must be designed carefully in order to
avoid nonpositive predictions and to accurately account for the singular
behavior in regions where the solution becomes small.

To address these difficulties, we develop a Physics-Informed Neural Network
(PINN) formulation in which the solution is represented by a parameterized
neural network and the governing differential equation is enforced through
its residual at a set of interior collocation points. The required spatial
derivatives of the neural approximation are evaluated using automatic
differentiation, so that the differential operator can be incorporated
directly into the optimization problem without requiring labeled solution
data.

In contrast with a standard PINN formulation based solely on residual
minimization and boundary penalties, the proposed framework is specifically
adapted to the structure of the singular elliptic problem. It combines three
main components:

\begin{itemize}

    \item a hard-constraint formulation that incorporates the Dirichlet
    boundary conditions directly into the neural representation;

    \item a positivity-preserving transformation that guarantees strictly
    positive predictions in the interior of the computational domain,
    ensuring that the singular term remains well-defined during training;

    \item a singularity-aware weighted residual loss that increases the
    contribution of regions where the predicted solution becomes small and
    the nonlinear singularity becomes more pronounced.

\end{itemize}

More precisely, let $u_{\theta}:\Omega\rightarrow\mathbb{R}$ denote the raw
output of a fully connected neural network with trainable parameters
$\theta$. Rather than using $u_{\theta}$ directly as the approximation of the
solution, we construct a constrained representation of the form

\begin{equation}
\hat{u}(x;\theta)
=
\eta(x)\,
\Phi\!\left(u_{\theta}(x)\right),
\label{eq:general_constrained_representation}
\end{equation}

where $\eta$ is a boundary-vanishing function satisfying

\begin{equation}
\eta(x)=0,
\qquad x\in\partial\Omega,
\label{eq:eta_boundary}
\end{equation}

and

\begin{equation}
\eta(x)>0,
\qquad x\in\Omega.
\label{eq:eta_positive}
\end{equation}

The mapping $\Phi:\mathbb{R}\rightarrow(0,\infty)$ is chosen to ensure
positivity of the neural approximation. In the present work, $\Phi$ is
implemented using the Softplus transformation. Consequently, the
representation~\eqref{eq:general_constrained_representation} satisfies both

\begin{equation}
\hat{u}(x;\theta)=0,
\qquad x\in\partial\Omega,
\end{equation}

and

\begin{equation}
\hat{u}(x;\theta)>0,
\qquad x\in\Omega,
\end{equation}

by construction.

For a given constrained approximation $\hat{u}(x;\theta)$, the differential
residual associated with problem~\eqref{eq:model} is

\begin{equation}
R(x;\theta)
=
-\Delta\hat{u}(x;\theta)
-
\frac{1}{\hat{u}(x;\theta)^{\alpha}}.
\label{eq:methodology_residual}
\end{equation}

The network parameters are determined by minimizing an objective function
constructed from evaluations of this residual at interior collocation points.
Because the magnitude of the singular nonlinear term depends strongly on the
predicted solution itself, a standard mean-squared residual may not assign
sufficient emphasis to regions where $\hat{u}$ becomes small. This
observation motivates the singularity-aware weighting strategy developed in
the following subsections.

The complete methodology therefore follows the sequence

\[
\boxed{
\begin{aligned}
u_{\theta}
&\;\longrightarrow\;
\text{positivity transformation}
\;\longrightarrow\;
\text{hard boundary constraint}
\;\longrightarrow\;
\hat{u}
\\[2mm]
&\;\longrightarrow\;
R
\;\longrightarrow\;
\text{weighted residual optimization}.
\end{aligned}
}
\]

The following subsections describe the neural architecture, hard enforcement
of the boundary conditions, positivity-preserving construction, residual
weighting strategy, and optimization procedure in detail.

\subsection{Neural Network Approximation}
\label{subsec:nn_approximation}

Following the standard PINN framework, we represent the unconstrained neural
output by a fully connected feed-forward neural network

\begin{equation}
u_{\theta}(x)
=
W^{L+1}z^{L}+b^{L+1},
\label{eq:nn_output}
\end{equation}

where the hidden representations are recursively defined by

\begin{equation}
z^{0}=x,
\end{equation}

and

\begin{equation}
z^{l}
=
\sigma
\left(
W^{l}z^{l-1}+b^{l}
\right),
\qquad
l=1,\ldots,L.
\label{eq:hidden_layers}
\end{equation}

Here, $L$ denotes the number of hidden layers,
$W^{l}$ and $b^{l}$ are the weight matrix and bias vector associated with
the $l$-th layer, respectively, and the complete set of trainable parameters
is denoted by

\begin{equation}
\theta
=
\left\{
W^{l},b^{l}
\right\}_{l=1}^{L+1}.
\label{eq:network_parameters}
\end{equation}

The hyperbolic tangent function is employed as the activation function in
the hidden layers,

\begin{equation}
\sigma(z)=\tanh(z).
\label{eq:tanh_activation}
\end{equation}

The use of a smooth activation function is particularly important in the
present setting because the elliptic operator involves second-order spatial
derivatives. The derivatives of the neural approximation required for the
evaluation of the Laplacian are computed using automatic differentiation.

It is important to distinguish the raw neural network output
$u_{\theta}(x)$ from the final approximation $\hat{u}(x;\theta)$ used in the
differential residual. The unconstrained output $u_{\theta}$ does not, by
itself, satisfy the boundary condition or guarantee positivity. These
properties are incorporated through the constrained representation

\begin{equation}
\hat{u}(x;\theta)
=
\eta(x)\,
\mathrm{Softplus}
\left(
u_{\theta}(x)
\right),
\label{eq:constrained_nn_recall}
\end{equation}

as described in the following subsections.

Consequently, the role of the fully connected neural network is to learn an
unconstrained latent representation, while the transformation applied to its
output ensures that the final approximation is compatible with the
structural requirements of the weakly singular elliptic problem.
\subsection{Hard Boundary Constraint}
\label{subsec:hard_boundary}

The weakly singular elliptic problem considered in this work is subject to
homogeneous Dirichlet boundary conditions. Rather than enforcing these
conditions through an additional penalty term in the loss function, we
incorporate them directly into the neural representation. This hard-constraint
strategy ensures that the boundary conditions are satisfied exactly for any
value of the trainable parameters.

Let $\eta:\overline{\Omega}\rightarrow\mathbb{R}$ be an auxiliary function
satisfying

\begin{equation}
\eta(x)=0,
\qquad x\in\partial\Omega,
\label{eq:eta_zero}
\end{equation}

and

\begin{equation}
\eta(x)>0,
\qquad x\in\Omega.
\label{eq:eta_interior}
\end{equation}

The final neural approximation is constructed as

\begin{equation}
\hat{u}(x;\theta)
=
\eta(x)\,
\mathrm{Softplus}
\left(
u_{\theta}(x)
\right),
\label{eq:hard_constraint}
\end{equation}

where $u_{\theta}(x)$ denotes the unconstrained output of the neural network
introduced in Section~\ref{subsec:nn_approximation}, and

\begin{equation}
\mathrm{Softplus}(z)
=
\ln\left(1+e^{z}\right).
\label{eq:softplus}
\end{equation}

By construction, for every $x\in\partial\Omega$,

\begin{equation}
\hat{u}(x;\theta)
=
\eta(x)\,
\mathrm{Softplus}\left(u_{\theta}(x)\right)
=
0,
\end{equation}

independently of the network parameters $\theta$. Hence, no separate
boundary-loss term is required for the homogeneous Dirichlet condition.

For the one-dimensional domain

\begin{equation}
\Omega=(0,1),
\end{equation}

a natural choice is

\begin{equation}
\eta(x)=x(1-x),
\label{eq:eta_1d}
\end{equation}

which is positive for $x\in(0,1)$ and vanishes exactly at $x=0$ and $x=1$.
The resulting approximation therefore takes the explicit form

\begin{equation}
\hat{u}(x;\theta)
=
x(1-x)\,
\mathrm{Softplus}
\left(
u_{\theta}(x)
\right).
\label{eq:hard_constraint_1d}
\end{equation}

For multidimensional domains, the function $\eta$ can be chosen according
to the geometry of $\Omega$, provided that it vanishes on
$\partial\Omega$ and remains positive in the interior.

In addition to enforcing the boundary condition, the representation
\eqref{eq:hard_constraint} preserves positivity in the interior. Indeed,
the Softplus function satisfies

\begin{equation}
\mathrm{Softplus}(z)>0,
\qquad z\in\mathbb{R},
\end{equation}

and therefore, since $\eta(x)>0$ for $x\in\Omega$,

\begin{equation}
\hat{u}(x;\theta)>0,
\qquad x\in\Omega.
\label{eq:positive_approximation}
\end{equation}

This property is essential for the singular elliptic problem because it
ensures that

\begin{equation}
\frac{1}{\hat{u}(x;\theta)^{\alpha}}
\end{equation}

remains well-defined at every interior collocation point throughout the
optimization process.

The hard-constrained representation thus simultaneously satisfies the
homogeneous Dirichlet boundary condition and preserves the positivity
required by the singular nonlinear source term.

\subsection{Physics-Informed Residual}
\label{subsec:physics_residual}

Once the constrained neural approximation $\hat{u}(x;\theta)$ has been
constructed, the governing differential equation is incorporated into the
learning process through its residual.

For a point
\[
x=(x_1,\ldots,x_d)\in\Omega,
\]
the Laplacian of the neural approximation is given by

\begin{equation}
\Delta\hat{u}(x;\theta)
=
\sum_{i=1}^{d}
\frac{\partial^{2}\hat{u}(x;\theta)}
{\partial x_i^{2}}.
\label{eq:nn_laplacian}
\end{equation}

The required spatial derivatives are evaluated using automatic
differentiation applied to the complete constrained approximation
$\hat{u}(x;\theta)$. This allows the derivatives of both the neural network
output and the transformations used to enforce the structural constraints
to be taken into account consistently.

The differential residual associated with problem~\eqref{eq:model} is then
defined as

\begin{equation}
R(x;\theta)
=
-\Delta\hat{u}(x;\theta)
-
\frac{1}
{\hat{u}(x;\theta)^{\alpha}}.
\label{eq:physics_residual}
\end{equation}

If $\hat{u}$ exactly satisfies the governing equation, then

\begin{equation}
R(x;\theta)=0,
\qquad x\in\Omega.
\end{equation}

In practice, the residual is evaluated at a set of interior collocation
points

\begin{equation}
\mathcal{X}_{r}
=
\left\{
x_i
\right\}_{i=1}^{N_r}
\subset\Omega,
\label{eq:collocation_set}
\end{equation}

and the trainable parameters $\theta$ are optimized so as to reduce the
residual over these points.

A standard PINN formulation would minimize the mean-squared residual

\begin{equation}
\mathcal{L}_{\mathrm{std}}(\theta)
=
\frac{1}{N_r}
\sum_{i=1}^{N_r}
R(x_i;\theta)^2.
\label{eq:standard_residual_loss}
\end{equation}

For the weakly singular problem considered here, however, the nonlinear
contribution

\begin{equation}
\frac{1}{\hat{u}(x;\theta)^{\alpha}}
\end{equation}

increases as the predicted solution becomes small. Consequently, different
regions of the computational domain may exhibit substantially different
levels of sensitivity to the singular nonlinearity.

This observation motivates the introduction of a singularity-aware weighted
residual loss, described in the following subsection, in order to place
additional emphasis on regions where $\hat{u}$ approaches zero.

\subsection{Weighted Physics-Informed Loss Function}
\label{subsec:weighted_loss}

The trainable parameters of the neural network are determined by minimizing
a residual-based objective function over the set of interior collocation
points
\[
\mathcal{X}_r
=
\left\{
x_i
\right\}_{i=1}^{N_r}
\subset\Omega.
\]

For a standard PINN formulation, all collocation points contribute equally
to the mean-squared residual. However, in the weakly singular problem
considered here, the magnitude of the nonlinear term

\begin{equation}
\frac{1}{\hat{u}(x;\theta)^{\alpha}}
\end{equation}

increases as the predicted solution $\hat{u}(x;\theta)$ approaches zero.
This motivates a residual weighting strategy that assigns additional
importance to regions where the singular nonlinear contribution is large.

We therefore define the singularity-aware weighted loss as

\begin{equation}
\mathcal{L}_{\mathrm{w}}(\theta)
=
\frac{1}{N_r}
\sum_{i=1}^{N_r}
w(x_i;\theta)
\,R(x_i;\theta)^2,
\label{eq:weighted_loss}
\end{equation}

where the weighting function is given by

\begin{equation}
w(x_i;\theta)
=
1+
\frac{\beta}
{\hat{u}(x_i;\theta)^{\alpha}},
\qquad
\beta>0.
\label{eq:weight_function}
\end{equation}

Here, $\beta$ is a positive weighting parameter controlling the strength of
the singularity-aware correction. When $\beta$ is small, the weighted
formulation remains close to the standard residual loss, whereas increasing
$\beta$ enhances the relative contribution of collocation points at which
the predicted solution is small.

The structure of~\eqref{eq:weight_function} is directly motivated by the
singular term appearing in the governing equation. Since

\begin{equation}
\hat{u}(x_i;\theta)\rightarrow0^{+}
\quad\Longrightarrow\quad
\hat{u}(x_i;\theta)^{-\alpha}\rightarrow+\infty,
\end{equation}

the weight increases automatically in regions where the predicted solution
approaches zero. Conversely, where $\hat{u}$ remains sufficiently away from
zero, the additional weighting contribution is less pronounced.

Moreover, because $\beta>0$ and $\hat{u}(x_i;\theta)>0$ at every interior
collocation point, we have

\begin{equation}
w(x_i;\theta)>1.
\label{eq:weight_lower_bound}
\end{equation}

Thus, no interior collocation point is removed from the residual
minimization; rather, the weighting modifies their relative contributions
according to the magnitude of the predicted solution.

For comparison, the corresponding standard PINN objective is

\begin{equation}
\mathcal{L}_{\mathrm{std}}(\theta)
=
\frac{1}{N_r}
\sum_{i=1}^{N_r}
R(x_i;\theta)^2.
\label{eq:standard_loss}
\end{equation}

The proposed formulation can therefore be interpreted as a
solution-dependent modification of the standard residual objective:

\begin{equation}
\mathcal{L}_{\mathrm{w}}(\theta)
=
\frac{1}{N_r}
\sum_{i=1}^{N_r}
\left(
1+
\frac{\beta}
{\hat{u}(x_i;\theta)^{\alpha}}
\right)
R(x_i;\theta)^2.
\label{eq:expanded_weighted_loss}
\end{equation}

The network parameters are finally determined through the optimization
problem

\begin{equation}
\theta^{*}
=
\operatorname*{arg\,min}_{\theta}
\mathcal{L}_{\mathrm{w}}(\theta).
\label{eq:optimization_problem}
\end{equation}

The purpose of this weighting strategy is to increase the influence, during
optimization, of regions where the singular nonlinear term is more
pronounced. Its effectiveness relative to the standard residual formulation
is evaluated numerically in Section~4.
\subsection{Training Procedure}
\label{subsec:training_procedure}

The overall training procedure of the proposed framework is summarized in
Algorithm~\ref{alg:pinn}. The method begins by generating a set of interior
collocation points and initializing the trainable parameters of the neural
network. The hard boundary construction and the positivity-preserving
Softplus transformation are then incorporated into the neural representation,
so that the homogeneous Dirichlet boundary condition is satisfied exactly
while the predicted solution remains strictly positive in the interior of
the computational domain.

At each training iteration, the Laplacian of the constrained approximation
is evaluated using automatic differentiation. The differential residual is
then computed at the interior collocation points and used to construct the
singularity-aware weighted loss introduced in
Section~\ref{subsec:weighted_loss}.

The optimization is performed in two successive stages. The Adam optimizer
is first employed to explore the parameter space and obtain a suitable
approximation. The resulting parameters are subsequently refined using the
L-BFGS optimizer, which is well suited to the final deterministic
optimization stage.

\begin{algorithm}[H]
\caption{Training procedure for the proposed PINN framework}
\label{alg:pinn}
\begin{algorithmic}[1]

\State Generate interior collocation points
$\mathcal{X}_r=\{x_i\}_{i=1}^{N_r}\subset\Omega$.

\State Initialize the neural network parameters $\theta$.

\State Construct the constrained neural approximation
\[
\hat{u}(x;\theta)
=
\eta(x)\,
\mathrm{Softplus}
\left(
u_{\theta}(x)
\right).
\]

\For{$k=1,\ldots,N_{\mathrm{Adam}}$}

    \State Evaluate $\hat{u}(x_i;\theta)$ at the collocation points.

    \State Compute the Laplacian
    $\Delta\hat{u}(x_i;\theta)$
    using automatic differentiation.

    \State Evaluate the residual
    \[
    R(x_i;\theta)
    =
    -\Delta\hat{u}(x_i;\theta)
    -
    \frac{1}
    {\hat{u}(x_i;\theta)^{\alpha}}.
    \]

    \State Compute the weights
    \[
    w(x_i;\theta)
    =
    1+
    \frac{\beta}
    {\hat{u}(x_i;\theta)^{\alpha}}.
    \]

    \State Compute the weighted loss
    \[
    \mathcal{L}_{\mathrm{w}}(\theta)
    =
    \frac{1}{N_r}
    \sum_{i=1}^{N_r}
    w(x_i;\theta)
    R(x_i;\theta)^2.
    \]

    \State Update $\theta$ using the Adam optimizer.

\EndFor

\State Use the parameters obtained from Adam as the initialization for
L-BFGS.

\State Refine $\theta$ by minimizing
$\mathcal{L}_{\mathrm{w}}(\theta)$
using L-BFGS.

\State \Return trained parameters $\theta^{*}$.

\end{algorithmic}
\end{algorithm}


\begin{figure}[ht]
\centering

\includegraphics[width=0.72\textwidth]{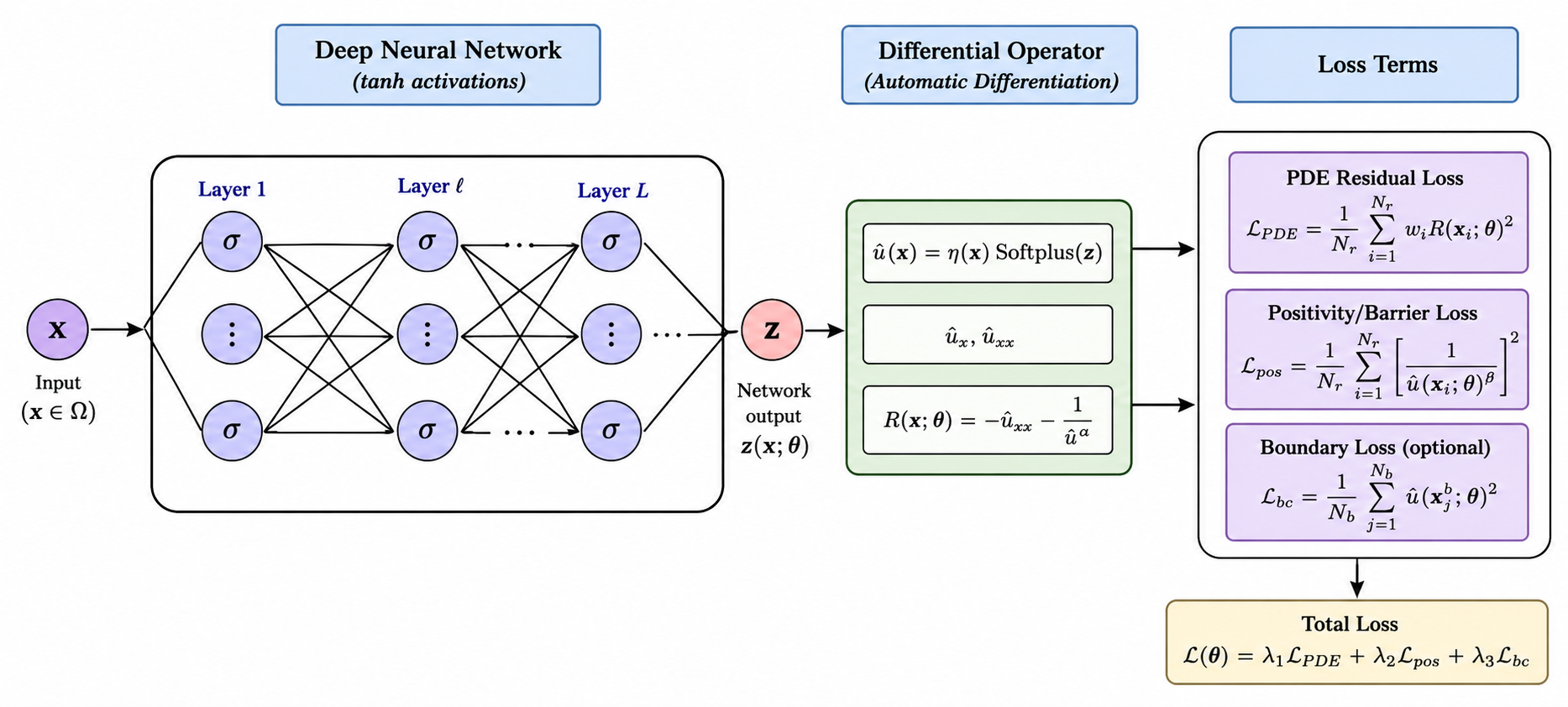}

\caption{
Workflow of the proposed PINN framework for weakly singular semilinear
elliptic equations. Interior collocation points are passed through the neural
network, followed by the positivity-preserving transformation and the hard
boundary construction. Automatic differentiation is then used to evaluate
the differential operator, and the resulting residual is incorporated into
the singularity-aware weighted loss used for optimization.
}

\label{fig:workflow}
\end{figure}

Figure~\ref{fig:workflow} summarizes the principal components of the proposed
methodology. Starting from the interior collocation points, the fully
connected neural network produces the unconstrained output
$u_{\theta}(x)$. The Softplus transformation and the boundary-vanishing
function $\eta(x)$ are subsequently applied to construct the constrained
approximation

\[
\hat{u}(x;\theta)
=
\eta(x)
\mathrm{Softplus}
\left(
u_{\theta}(x)
\right).
\]

Automatic differentiation is then employed to evaluate the Laplacian of the
complete constrained approximation. This leads to the differential residual

\[
R(x;\theta)
=
-\Delta\hat{u}(x;\theta)
-
\frac{1}
{\hat{u}(x;\theta)^{\alpha}},
\]

which is incorporated into the weighted objective function. The network
parameters are first optimized using Adam and subsequently refined using
L-BFGS.

After training, the objective is to obtain parameters $\theta^{*}$ for which

\[
R(x;\theta^{*})\approx0
\]

throughout the computational domain, while the architectural construction
guarantees

\[
\hat{u}(x;\theta^{*})>0,
\qquad x\in\Omega,
\]

and

\[
\hat{u}(x;\theta^{*})=0,
\qquad x\in\partial\Omega.
\]

The positivity property ensures that the singular nonlinear contribution

\[
\frac{1}
{\hat{u}(x;\theta^{*})^{\alpha}}
\]

remains well-defined at interior collocation points.


\begin{remark}

For comparison, a standard residual-based PINN minimizes

\begin{equation}
\mathcal{L}_{\mathrm{std}}(\theta)
=
\frac{1}{N_r}
\sum_{i=1}^{N_r}
R(x_i;\theta)^2.
\label{eq:standard_loss_training}
\end{equation}

In this formulation, all collocation points enter the objective function
without an explicit solution-dependent weighting factor. For the weakly
singular equation considered here, however,

\[
\hat{u}(x;\theta)\rightarrow0^{+}
\]

implies

\[
\hat{u}(x;\theta)^{-\alpha}
\rightarrow+\infty.
\]

Consequently, regions in which the predicted solution becomes small are
particularly sensitive to the singular nonlinear contribution.

The proposed weighted formulation modifies the standard objective according
to

\begin{equation}
\mathcal{L}_{\mathrm{w}}(\theta)
=
\frac{1}{N_r}
\sum_{i=1}^{N_r}
\left(
1+
\frac{\beta}
{\hat{u}(x_i;\theta)^{\alpha}}
\right)
R(x_i;\theta)^2,
\qquad
\beta>0.
\label{eq:weighted_loss_training}
\end{equation}

The additional weighting factor therefore increases the relative
contribution of collocation points for which the predicted solution is
small. The purpose of this construction is to direct additional
optimization attention toward regions where the singular nonlinear term is
more pronounced.

Figure~\ref{fig:loss_compare} provides a schematic comparison between the
standard residual formulation and the proposed weighted formulation. At this
stage, the curves are included only to illustrate the type of convergence
comparison used in the numerical study. They must therefore be replaced by
the actual training histories obtained from the experiments reported in
Section~4 before the final version of the manuscript is prepared.


\begin{figure}[H]
\centering

\begin{tikzpicture}

\begin{axis}[
width=0.78\textwidth,
height=6cm,
xlabel={Iterations ($\times10^4$)},
ylabel={Illustrative Loss (dB)},
xmin=0,
xmax=15,
ymin=-35,
ymax=0,
grid=both,
grid style={dashed,gray!30},
legend style={
    at={(0.98,0.98)},
    anchor=north east,
    font=\small
},
tick label style={font=\small},
label style={font=\small}
]

\addplot[
blue,
thick
]
coordinates{
(0,0)
(1,-5)
(2,-11)
(3,-17)
(4,-21)
(5,-24)
(6,-26)
(7,-27)
(8,-28)
(9,-29)
(10,-30)
(11,-31)
(12,-31.5)
(13,-32)
(14,-32.3)
(15,-32.5)
};

\addlegendentry{Weighted PINN}

\addplot[
red,
dashed,
thick
]
coordinates{
(0,0)
(1,-3)
(2,-6)
(3,-8)
(4,-9.5)
(5,-10.5)
(6,-11)
(7,-11.2)
(8,-11.4)
(9,-11.5)
(10,-11.6)
(11,-11.7)
(12,-11.8)
(13,-11.9)
(14,-12)
(15,-12.1)
};

\addlegendentry{Standard PINN}

\end{axis}

\end{tikzpicture}

\caption{
Schematic comparison of the convergence behavior of the standard residual
formulation and the proposed singularity-aware weighted formulation. The
curves are illustrative placeholders and are not numerical results. They
will be replaced by the actual convergence histories obtained from the
experiments reported in Section~4.
}

\label{fig:loss_compare}

\end{figure}
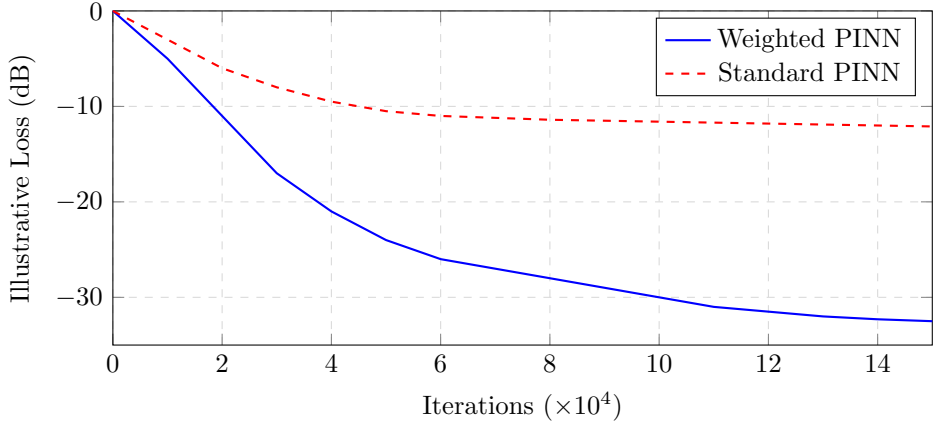

\end{remark}

\subsection{Evaluation Metrics}
\label{subsec:evaluation_metrics}

The performance of the proposed PINN framework is evaluated using several
quantitative metrics designed to assess both solution accuracy and
satisfaction of the governing differential equation.

Whenever an analytical solution is available, the neural approximation is
compared directly with the exact solution. For problems for which no explicit
solution is available, a high-resolution numerical solution computed using a
classical finite element method is employed as a reference solution. In the
following, $u_{\mathrm{ref}}$ denotes either the exact solution or the
high-resolution numerical reference solution, depending on the considered
experiment.

The relative $L^2$ error is defined as

\begin{equation}
E_{L^2}
=
\frac{
\left\|
u_{\mathrm{ref}}-\hat{u}
\right\|_{L^2(\Omega)}
}{
\left\|
u_{\mathrm{ref}}
\right\|_{L^2(\Omega)}
}.
\label{eq:relative_l2_error}
\end{equation}

This metric provides a global measure of the discrepancy between the neural
approximation and the reference solution over the computational domain.

To quantify the largest pointwise discrepancy, we also consider the relative
$L^\infty$ error

\begin{equation}
E_{L^\infty}
=
\frac{
\left\|
u_{\mathrm{ref}}-\hat{u}
\right\|_{L^\infty(\Omega)}
}{
\left\|
u_{\mathrm{ref}}
\right\|_{L^\infty(\Omega)}
},
\label{eq:relative_linf_error}
\end{equation}

where

\begin{equation}
\left\|
u_{\mathrm{ref}}-\hat{u}
\right\|_{L^\infty(\Omega)}
=
\sup_{x\in\Omega}
\left|
u_{\mathrm{ref}}(x)-\hat{u}(x)
\right|.
\end{equation}

In the numerical implementation, this quantity is approximated over a
sufficiently dense set of evaluation points.

In addition to solution-based errors, the satisfaction of the governing
equation is assessed through the mean-squared residual. Let

\begin{equation}
\mathcal{X}_{\mathrm{test}}
=
\left\{
x_j^{\mathrm{test}}
\right\}_{j=1}^{N_{\mathrm{test}}}
\subset\Omega
\end{equation}

denote a set of evaluation points independent of the collocation points used
during training. The residual error is defined by

\begin{equation}
\mathrm{MSE}_{R}
=
\frac{1}{N_{\mathrm{test}}}
\sum_{j=1}^{N_{\mathrm{test}}}
R\left(x_j^{\mathrm{test}};\theta^{*}\right)^2,
\label{eq:residual_mse}
\end{equation}

where $\theta^{*}$ denotes the optimized network parameters.

Evaluating the residual on an independent set of points provides a
complementary measure of how well the trained neural approximation satisfies
the differential equation beyond the collocation points explicitly used in
the optimization process.

Finally, the convergence behavior of the optimization procedure is analyzed
through the evolution of the training objectives. In particular, the
weighted loss

\begin{equation}
\mathcal{L}_{\mathrm{w}}(\theta)
=
\frac{1}{N_r}
\sum_{i=1}^{N_r}
\left(
1+
\frac{\beta}
{\hat{u}(x_i;\theta)^{\alpha}}
\right)
R(x_i;\theta)^2
\end{equation}

is monitored during training. When comparing the proposed method with the
standard PINN formulation, the corresponding standard residual loss
$\mathcal{L}_{\mathrm{std}}$ is monitored under the same experimental
conditions.

The combination of the relative $L^2$ error, relative $L^\infty$ error,
independent residual error, and optimization histories provides a
comprehensive assessment of the accuracy and convergence behavior of the
proposed singularity-aware PINN framework.

\section{Numerical Experiments}
\label{sec:numerical_experiments}

In this section, we investigate the numerical performance of the proposed
PINN framework for the approximation of weakly singular semilinear elliptic
equations of the form

\begin{equation}
-\Delta u
=
\frac{1}{u^{\alpha}},
\qquad
0<\alpha<1,
\label{eq:numerical_model}
\end{equation}

subject to homogeneous Dirichlet boundary conditions and the interior
positivity requirement.

The numerical experiments are designed to assess the accuracy and
convergence behavior of the proposed singularity-aware weighted formulation.
Particular attention is given to regions where the solution becomes small,
since

\[
u(x)\rightarrow0^{+}
\quad\Longrightarrow\quad
u(x)^{-\alpha}\rightarrow+\infty.
\]

These regions therefore provide a particularly relevant test of the
positivity-preserving architecture and the proposed residual weighting
strategy.

Throughout the experiments, the solution is approximated using the fully
connected feed-forward neural network introduced in
Section~\ref{subsec:nn_approximation}. Unless otherwise specified, the hidden
layers employ the hyperbolic tangent activation function. The raw neural
network output $u_{\theta}$ is transformed according to

\begin{equation}
\hat{u}(x;\theta)
=
\eta(x)\,
\mathrm{Softplus}
\left(
u_{\theta}(x)
\right),
\label{eq:numerical_constrained_solution}
\end{equation}

so that the homogeneous Dirichlet boundary condition is satisfied exactly
and the approximation remains strictly positive in the interior of the
domain.

The network parameters are optimized using the two-stage training strategy
described in Section~\ref{subsec:training_procedure}. Adam is first employed
during the initial optimization stage, after which L-BFGS is used to refine
the resulting approximation.

For the principal one-dimensional experiment, the baseline configuration is
summarized in Table~\ref{tab:hyper}. Additional experiments may modify
individual parameters in order to study their influence on the numerical
performance.

\begin{table}[H]
\centering
\caption{Baseline configuration used for the one-dimensional weakly singular
problem.}
\label{tab:hyper}

\begin{tabular}{lc}
\toprule
Parameter & Value \\
\midrule
Singularity exponent $\alpha$ & $0.5$ \\
Interior collocation points $N_r$ & $5000$ \\
Hidden layers & $3$ \\
Neurons per hidden layer & $24$ \\
Hidden activation & $\tanh$ \\
Output transformation & Softplus \\
Boundary function & $x(1-x)$ \\
Optimization & Adam followed by L-BFGS \\
\bottomrule
\end{tabular}

\end{table}

The same baseline configuration is used when comparing the standard
residual loss with the proposed weighted residual loss. In this way, the
effect of the weighting strategy can be assessed without simultaneously
changing the neural architecture or the number of collocation points.

Whenever an analytical solution is available, the neural approximation is
compared directly with the exact solution. When no explicit analytical
solution is available, a sufficiently accurate numerical solution computed
using a classical numerical method is employed as the reference solution.
For the higher-dimensional experiments, a high-resolution finite element
solution is used for this purpose.

The numerical performance is assessed using the metrics introduced in
Section~\ref{subsec:evaluation_metrics}, including the relative $L^2$ error,
the relative $L^\infty$ error, and the mean-squared differential residual
evaluated on an independent set of test points. Training histories are also
examined to compare the convergence behavior of the standard and weighted
formulations.

The experiments are organized to address three main questions:

\begin{enumerate}

    \item Can the hard-constrained, positivity-preserving neural
    representation accurately approximate the solution of the weakly
    singular elliptic problem?

    \item Does the proposed singularity-aware weighting improve the
    numerical performance relative to the standard residual formulation
    under the same network architecture and training conditions?

    \item How does the behavior of the method change as the singularity
    exponent $\alpha$ and the dimension of the computational domain are
    varied?

\end{enumerate}

The following subsections first consider the one-dimensional weakly singular
problem, followed by a controlled manufactured-solution verification test.
Additional experiments investigate the influence of the singularity
parameter and the extension of the method to higher-dimensional domains.

\subsection{One-Dimensional Weakly Singular Elliptic Problem}
\label{subsec:1d_singular}

We first investigate the proposed neural framework on a one-dimensional
weakly singular semilinear elliptic problem. Let $\Omega=(0,1)$ and consider

\begin{equation}
-\frac{d^{2}u}{dx^{2}}
=
\frac{1}{u^{\alpha}},
\qquad x\in(0,1),
\label{eq:example1}
\end{equation}

subject to the homogeneous Dirichlet boundary conditions

\begin{equation}
u(0)=u(1)=0,
\label{eq:bc1}
\end{equation}

together with the positivity condition

\begin{equation}
u(x)>0,
\qquad x\in(0,1),
\label{eq:positivity1}
\end{equation}

where

\begin{equation}
0<\alpha<1.
\label{eq:alpha_range}
\end{equation}

This range corresponds to the weakly singular regime considered throughout
this work. Although the solution remains strictly positive in the interior
of the domain, the homogeneous boundary conditions imply that
$u(x)\rightarrow0^{+}$ as $x$ approaches either endpoint. Consequently,

\begin{equation}
u(x)^{-\alpha}\rightarrow+\infty,
\end{equation}

and the nonlinear source term becomes unbounded near the boundary. This
behavior constitutes the main numerical difficulty of the problem and
provides a natural benchmark for evaluating the proposed method.

For the first numerical experiment, we set

\begin{equation}
\alpha=0.5.
\label{eq:alpha_example1}
\end{equation}

To incorporate the boundary conditions and positivity requirement directly
into the approximation, we employ the constrained representation introduced
in Section~\ref{subsec:hard_boundary}. In one dimension, the boundary
function is chosen as

\begin{equation}
\eta(x)=x(1-x),
\end{equation}

which leads to the neural approximation

\begin{equation}
\hat{u}(x;\theta)
=
x(1-x)\,
\mathrm{Softplus}
\left(
u_{\theta}(x)
\right).
\label{eq:hard_1d}
\end{equation}

Since

\[
x(1-x)=0
\qquad
\text{for }x\in\{0,1\},
\]

the homogeneous Dirichlet boundary conditions are satisfied exactly:

\begin{equation}
\hat{u}(0;\theta)
=
\hat{u}(1;\theta)
=
0.
\end{equation}

Moreover, since $x(1-x)>0$ for $x\in(0,1)$ and the Softplus transformation
is strictly positive,

\begin{equation}
\hat{u}(x;\theta)>0,
\qquad x\in(0,1).
\end{equation}

Thus, both the boundary conditions and the interior positivity requirement
are embedded directly into the neural representation without introducing
additional boundary penalty terms into the objective function.

For the approximation~\eqref{eq:hard_1d}, the differential residual is
defined by

\begin{equation}
R(x;\theta)
=
-\frac{d^{2}\hat{u}(x;\theta)}{dx^{2}}
-
\frac{1}{\hat{u}(x;\theta)^{\alpha}}.
\label{eq:residual_1d}
\end{equation}

The required spatial derivatives are evaluated using automatic
differentiation.

The proposed formulation minimizes the singularity-aware weighted loss

\begin{equation}
\mathcal{L}_{\mathrm{w}}(\theta)
=
\frac{1}{N_r}
\sum_{i=1}^{N_r}
\left(
1+
\frac{\beta}
{\hat{u}(x_i;\theta)^{\alpha}}
\right)
R(x_i;\theta)^2,
\label{eq:weighted_loss_1d}
\end{equation}

where $\beta>0$ controls the strength of the residual weighting. The
additional factor increases the relative contribution of collocation points
for which the predicted solution is small and the singular nonlinear
contribution is therefore more pronounced.

For comparison, we also consider the standard residual objective

\begin{equation}
\mathcal{L}_{\mathrm{std}}(\theta)
=
\frac{1}{N_r}
\sum_{i=1}^{N_r}
R(x_i;\theta)^2.
\label{eq:standard_loss_1d}
\end{equation}

This comparison is designed to assess the contribution of the proposed
singularity-aware weighting independently of the neural architecture and
hard-constrained representation.

The neural network consists of three hidden layers with twenty-four neurons
per layer and hyperbolic tangent activation functions. A total of
$N_r=5000$ interior collocation points are employed during training. The
network parameters are optimized using the two-stage procedure described in
Section~\ref{subsec:training_procedure}: Adam is first employed during the
initial optimization stage, after which L-BFGS is used to refine the
approximation.

The same neural architecture, collocation points, initialization strategy,
and optimization procedure are used for the standard and weighted
formulations so that their numerical performance can be compared under
consistent computational conditions.

Since problem~\eqref{eq:example1} is not evaluated here using a prescribed
closed-form solution, an independently computed high-accuracy numerical
solution is employed as the reference solution. This reference solution is
used exclusively for post-training validation and does not participate in
the optimization of the neural network.

The approximation quality is assessed using the relative $L^2$ error

\begin{equation}
\mathcal{E}_{L^2}
=
\frac{
\left\|
u_{\mathrm{ref}}-\hat{u}
\right\|_{L^2(\Omega)}
}{
\left\|
u_{\mathrm{ref}}
\right\|_{L^2(\Omega)}
},
\label{eq:l2_error_1d}
\end{equation}

and the relative $L^\infty$ error

\begin{equation}
\mathcal{E}_{L^\infty}
=
\frac{
\left\|
u_{\mathrm{ref}}-\hat{u}
\right\|_{L^\infty(\Omega)}
}{
\left\|
u_{\mathrm{ref}}
\right\|_{L^\infty(\Omega)}
}.
\label{eq:linf_error_1d}
\end{equation}

In the numerical implementation, these quantities are evaluated on a
sufficiently dense set of test points. In addition, the mean-squared
differential residual is evaluated on an independent test set, and the
optimization histories of the standard and weighted formulations are
recorded.

This one-dimensional benchmark therefore serves two complementary purposes.
First, it evaluates the ability of the hard-constrained,
positivity-preserving neural representation to approximate a solution in
the presence of an unbounded nonlinear source term near the boundary.
Second, the comparison between~\eqref{eq:weighted_loss_1d} and
\eqref{eq:standard_loss_1d} isolates the effect of the proposed
singularity-aware residual weighting under identical computational
conditions.


\subsection{Manufactured-Solution Verification Test}
\label{subsec:manufactured}

The previous experiment addresses the original weakly singular equation.
However, an additional benchmark with an explicitly known solution is useful
for quantitatively verifying the numerical implementation. For this purpose,
we introduce a manufactured-solution problem for which the exact solution is
available in closed form.

We prescribe

\begin{equation}
u_{\mathrm{ex}}(x)
=
1+x(1-x),
\qquad x\in(0,1).
\label{eq:manufactured_solution}
\end{equation}

Its first and second derivatives are

\begin{equation}
u_{\mathrm{ex}}'(x)=1-2x,
\qquad
u_{\mathrm{ex}}''(x)=-2,
\end{equation}

and therefore

\begin{equation}
-u_{\mathrm{ex}}''(x)=2.
\label{eq:manufactured_second_derivative}
\end{equation}

We then consider the modified semilinear problem

\begin{equation}
-\frac{d^{2}u}{dx^{2}}
=
\frac{1}{u^{\alpha}}
+
f(x),
\qquad x\in(0,1),
\label{eq:manufactured_problem}
\end{equation}

where the forcing term is defined by

\begin{equation}
f(x)
=
2-
\frac{1}
{\left[1+x(1-x)\right]^{\alpha}}.
\label{eq:manufactured_source}
\end{equation}

Substituting~\eqref{eq:manufactured_solution} into
\eqref{eq:manufactured_problem} gives

\begin{align}
\frac{1}{u_{\mathrm{ex}}(x)^{\alpha}}
+
f(x)
&=
\frac{1}
{\left[1+x(1-x)\right]^{\alpha}}
+
2
-
\frac{1}
{\left[1+x(1-x)\right]^{\alpha}}
\nonumber\\
&=
2
=
-u_{\mathrm{ex}}''(x),
\end{align}

which confirms that~\eqref{eq:manufactured_solution} is an exact solution of
the modified differential equation.

The corresponding Dirichlet boundary conditions are

\begin{equation}
u(0)=u(1)=1.
\label{eq:manufactured_bc}
\end{equation}

Because these boundary conditions are nonhomogeneous, the hard-constrained
representation used for the original singular problem must be modified. We
therefore define

\begin{equation}
\hat{u}(x;\theta)
=
1
+
x(1-x)\,
\mathrm{Softplus}
\left(
u_{\theta}(x)
\right).
\label{eq:manufactured_hard_constraint}
\end{equation}

Since $x(1-x)=0$ at $x=0$ and $x=1$, this representation satisfies

\begin{equation}
\hat{u}(0;\theta)
=
\hat{u}(1;\theta)
=
1
\end{equation}

exactly for every value of the network parameters. Moreover, since

\[
x(1-x)\geq0,
\qquad x\in[0,1],
\]

and the Softplus function is strictly positive, we have

\begin{equation}
\hat{u}(x;\theta)\geq1,
\qquad x\in[0,1].
\label{eq:manufactured_positive}
\end{equation}

The differential residual associated with the manufactured problem is

\begin{equation}
R_{\mathrm{m}}(x;\theta)
=
-\frac{d^{2}\hat{u}(x;\theta)}{dx^{2}}
-
\frac{1}{\hat{u}(x;\theta)^{\alpha}}
-
f(x).
\label{eq:manufactured_residual}
\end{equation}

For the exact solution,

\begin{equation}
u_{\mathrm{ex}}(x)
=
1+x(1-x)
\geq1,
\qquad x\in[0,1].
\end{equation}

Consequently,

\begin{equation}
0
<
u_{\mathrm{ex}}(x)^{-\alpha}
\leq1,
\qquad x\in[0,1],
\end{equation}

and the nonlinear term remains bounded throughout the computational domain.
The manufactured problem therefore does not reproduce the boundary-induced
singularity of the original problem considered in
Section~\ref{subsec:1d_singular}.

Its purpose is instead to provide a controlled verification test in which
the exact solution is known. This makes it possible to assess the accuracy
of the neural approximation independently of the numerical reference
solution required for the original singular problem.

Because $u_{\mathrm{ex}}$ is explicitly known, the approximation error can
be evaluated directly. The relative $L^2$ error is defined as

\begin{equation}
\mathcal{E}_{L^2}^{\mathrm{m}}
=
\frac{
\left\|
u_{\mathrm{ex}}-\hat{u}
\right\|_{L^2(\Omega)}
}{
\left\|
u_{\mathrm{ex}}
\right\|_{L^2(\Omega)}
},
\label{eq:manufactured_l2}
\end{equation}

while the relative $L^\infty$ error is given by

\begin{equation}
\mathcal{E}_{L^\infty}^{\mathrm{m}}
=
\frac{
\left\|
u_{\mathrm{ex}}-\hat{u}
\right\|_{L^\infty(\Omega)}
}{
\left\|
u_{\mathrm{ex}}
\right\|_{L^\infty(\Omega)}
}.
\label{eq:manufactured_linf}
\end{equation}

In the numerical implementation, the $L^\infty$ norm is approximated over a
sufficiently dense set of evaluation points. The mean-squared differential
residual is additionally evaluated on an independent test set in order to
measure how accurately the trained neural approximation satisfies the
manufactured differential equation away from the training collocation
points.

The manufactured-solution benchmark complements, rather than replaces, the
original singular experiment. The problem considered in
Section~\ref{subsec:1d_singular} evaluates the proposed method in the
presence of a genuinely unbounded nonlinear source term as the solution
approaches the boundary, whereas the present experiment provides a
controlled verification against an explicitly known analytical solution.
Together, the two experiments assess complementary aspects of the proposed neural solver.

\section*{Data Availability Statement}

The datasets generated or analyzed during the current study are available
from the corresponding author upon reasonable request.

\author{
Badr Oulgiht\\
Independent Researcher, France\\
\texttt{ba.oulgiht@gmail.com}
}
\end{document}